\documentclass[11pt,notitlepage,twoside]{article}
\usepackage{latexsym}
\usepackage{amsmath}
\usepackage{amsfonts}
\usepackage{amssymb}
\renewcommand{\a }{\alpha }

\renewcommand{\d}{\delta }

\newcommand{\D }{\Delta }

\newcommand{\e }{\varepsilon }

\renewcommand{\l }{\lambda }
\renewcommand{\L }{\Lambda }

\newcommand{\n }{\nabla }
\newcommand{\var }{\varphi }

\renewcommand{\th }{\theta }

\renewcommand{\O }{\Omega }

\newcommand{\be}{\begin{equation}}
\newcommand{\ee}{\end{equation}}
\newenvironment{pf}{\noindent{\bf Proof.}\enspace}{
\hfill$\Box$\medskip}
\newenvironment{pfn}[1]{\noindent{\bf Proof of {#1}\enspace}}{
\hfill$\Box$\medskip}
\newcommand{\R}{\mathbb{R}}
 
\newcommand{\N}{\mathbb{N}}

\newtheorem{thm}{Theorem}[section]
\newtheorem{pro}[thm]{Proposition}
\newtheorem{lem}[thm]{Lemma}
\newtheorem{rem}[thm]{Remark}
\newtheorem{cor}[thm]{Corollary}

\numberwithin{equation}{section}

\author{M. Ben Ayed$^a$, K. El Mehdi$^{b,c}$\thanks{Corresponding author. E-mails: \texttt{khalil@univ-nkc.mr} and \texttt{elmehdik@ictp.trieste.it}. }, M. Hammami$^a$  \& M. Ould Ahmedou$^d$ \thanks{E-mails: M. Ben Ayed: \texttt{Mohamed.Benayed@fss.rnu.tn},  M. Hammami: \texttt{Mokhless.Hammami@fss.rnu.tn}, M. Ould Ahmedou:  \texttt{ahmedou@everest.mathematik.uni-tuebingen.de}.} \\
{\footnotesize
a:  D\'epartement de Math\'ematiques, Facult{\'e} des Sciences de Sfax, Route Soukra, Sfax, Tunisia.}\\
{\footnotesize
b:  Facult\'e des Sciences et Techniques, Universit\'e de Nouakchott, Nouakchott, Mauritania.}\\
{\footnotesize
 c:  The Abdus Salam ICTP, Mathematics Section, Strada Costiera 11, 34014 Trieste, Italy.}\\
{\footnotesize
d: Mathematisches institut, Auf der Morgenstelle 10, D-72076 Tubingen, Germany.}
}

\title { \Large \textbf{On a  Yamabe Type
 Problem on Three Dimensional Thin Annulus  }}

\begin{document}

\date{ }

\maketitle
{\footnotesize
\noindent
{\bf ABSTRACT.-}
We consider the problem: $ (P_{\varepsilon}):\,
-\Delta u_\varepsilon =  u_\varepsilon^{5},\, u_\varepsilon
>0\,\,  \mbox{  in  }  A_\varepsilon; \,
u_\varepsilon= 0\,\, \mbox{ on }  \partial A_\varepsilon
$, where $\{A_\e \subset \R^3, \e >0\}$ is a family of bounded annulus shaped domains such that $A_\e$ becomes ``thin'' as $\e\to 0$. We show that, for
 any given constant $C>0,$  there exists $\varepsilon_0>0$
such that for any $\varepsilon <\varepsilon_0$, the problem $(P_\e)$
has no solution $u_\varepsilon,$  whose energy, $\int_{A_\varepsilon}|\nabla
u_\varepsilon |^2,$ is less than C. Such a result extends to dimension three
a result previously known in higher dimensions. Although the strategy to
prove this result is the same as in higher dimensions,  we need a more
careful and delicate blow up analysis of asymptotic profiles of solutions $u_\e$ when $\e\to 0$. 
}

\section{Introduction }
We consider  the following nonlinear elliptic problem
 $$
(P_\Omega )\quad \left\{
\begin{array}{ccccc}
 -\Delta u &=& u^{5}, & u >0 & \mbox{  in  }\,  \Omega\\
 u &=&0 & & \mbox{ on }\,  \partial  \Omega,
\end{array}
\right.
$$
where $\Omega$ is a smooth bounded  domain in $\R^3$.\\
The equation ($P_{\Omega}$) arises in many mathematical and physical contexts (see \cite{Br}), but its greatest
interest lies in its relation to the Yamabe problem. From this geometric
point of view, we think of $u$ as defining the conformal metric $g_{ij} =
u^{\frac{4}{n - 2}} \, \d_{ij}$. Equation ($P_{\O}$) then says that the metric
$g$ has constant scalar curvature.

It is well known that if $\Omega$ is starshaped, $(P_\Omega) $ has no
solution (see Pohozaev \cite{P}) and if $\Omega$ has nontrivial topology,
in the sense that $H_{2k-1}(\Omega ;Q) \neq 0$ or $H_{k}(\Omega ;Z/2Z) \neq
0$ for some $k\in\N,$ Bahri and Coron \cite{BC} have shown that  $(P_\Omega )
\,$ has a solution. Nevertheless, Ding \cite{Di} (see also Dancer \cite{Da})
gave the example of contractible domain on which $(P_\Omega)$  has a solution.
 Then, the question related to existence or nonexistence of solution of $(P_\O)$ remained open.

In this paper, we study the problem $(P_\Omega)$  when $\Omega  = A_{
\varepsilon}$ is an annulus-shaped domain in $\R^3$and $\varepsilon$ is a small positive parameter. More precisely, let $f$ be any smooth function:
$$
f:\R^{2}\longrightarrow [1,2]\, \, ,
(\theta_1 ,\theta_{2} ) \longrightarrow f(\theta_1, \theta_{2})
$$
\noindent
which is periodic of period $\pi$ with respect to  $\theta_1$
and of period $2\pi$ with respect to $\theta_{2}$.\\
We set
$$
S_1 (f) = \left\{ x\in\R^3 / r= f(\theta_1, \theta_2)  \right\},
$$
\noindent
where $ (r ,\theta_1, \theta_2)$ are the polar coordinates  of  $x$.\\
For $ \e $ positive small enough, we introduce the following map
$$
g_\e:S_1 (f)\longrightarrow g_\e (S_1 (f)) = S_2(f), \, \, \, x \longmapsto g_\e (x)= x + \e n_x,
$$
where $ n_x $ is the outward normal to $ S_1 (f) $ at $ x $.
We denote by $ (A_\e)_{\e > 0} $ the family of annulus-shaped domain in $ \R
^3 $ such that $ \partial A_\e = S_1 (f) \cup S_2 (f) $.

We are mainly interested in the existence of finite energy
solutions, our main result is the following Theorem.
\begin{thm}\label{t:11}
Let C be any positive constant. Then, there exists $\varepsilon_0
>0$ such that for any $\varepsilon <\varepsilon_0$,  the problem $
(P_{\varepsilon}):\, -\Delta u_\varepsilon =  u_\varepsilon^{5},\,
u_\varepsilon >0
 \mbox{  in  }  A_\varepsilon,\,\,
u_\varepsilon =0 \mbox{ on }  \partial A_\varepsilon
$,
has no solution such that $\int_{A_\varepsilon}|\nabla u_\varepsilon |^2 \leq C.$
\end{thm}

Such a nonexistence result of finite energy solutions to Yamabe
type problems on nontrivial domains is a new and interesting
phenomenon, and it is a subject of current investigations by the
authors. It turns out that such a nonexistence result of finite energy
solutions is closely related to nonexistence results of solutions of finite
Morse index, and has its explanation in the behavior of the first eigenvalue
of  Laplace operator, or more generally of Laplace Beltrami operator on complete
manifolds. We hope that such results will be useful to find  necessary and
sufficient conditions on the manifold for the solvability of Yamabe problem
on complete manifolds. The results of such investigations  will appear
elsewhere. We notice that the  higher dimensional analogue of our result has been
recently proved  by the first
three authors \cite{BEH}.\\
Our strategy to prove Theorem \ref{t:11} is the same as in higher
dimensions, however as usual in elliptic equations involving critical Sobolev
exponent \cite{BN}, we need  more refined estimates  of the
asymptotic profiles of solutions $u_\e$ when $\e\to 0$ to treat the three
dimensional case. Such
refined estimates, which are of self interest, are highly
nontrivial and uses in a crucial way the refined properties of
blowing up solutions of Yamabe type problems in the spirit of R. Schoen
\cite{Sch1}, \cite{Sch2}, \cite{S3} and
Y. Y. Li \cite{yyli}.
 The input of such a refined blow up analysis enables us to
rule out some {\it bad configurations} for which the higher
dimensional
estimates cannot be improved.\\
Another ingredient of our proof is a careful expansion of the
Euler Lagrange functional associated to $(P_\e)$, and its gradient
near a small neighborhood of highly concentrated functions. To perform such 
expansions we extensively make use of the techniques developed by A. Bahri \cite{B} and
O. Rey \cite{R}, \cite{R2} in the framework of the {\it Theory of
critical points at infinity}.

The organization of the paper is as follows. The next section is
devoted to set up  some notation. In Section 3, we study the
asymptotic behavior of bounded energy solutions of $(P_\e)$. In
section  4, we  prove Theorem \ref{t:11}. Lastly,  we prove in
Section 5 some useful facts and careful estimates needed for the
previous sections.
\section{Notation}
 We denote by $G_{\varepsilon}$ the Green's function of Laplace operator defined by
\begin{eqnarray}\label{e:21}
 \forall\, x\in A_{\varepsilon }\quad
-\Delta G_{\varepsilon}(x,.) = c'\delta_x \mbox{  in  }
A_{\varepsilon}\quad, \quad G_{\varepsilon}(x,.) = 0 \mbox{ on }
\partial  A_{\varepsilon},
\end{eqnarray}
 \noindent
where $\delta_x$ is the Dirac mass at $x$ and $c' = meas(S^{2}).$\\
We denote by  $H_{\varepsilon}$ the regular part of
$G_{\varepsilon},$  that is,
\begin{eqnarray}\label{e:22}
\quad H_{\varepsilon}(x_1 ,x_2) = {|x_1-x_2|^{-1}} - G_{\varepsilon}
(x_1 , x_2) ,\mbox{ for } (x_1 ,x_2)\in A_{\varepsilon}\times A_{\varepsilon}.
\end{eqnarray}
 \noindent
For $p\in\N^*$ and $ x= (x_1 ,...,x_p)\in A_{\varepsilon}^{p} $,
we denote by $M = M_{\varepsilon}(x)$ the matrix defined by
\begin{eqnarray}\label{e:23}
  M = (m_{ij})_{1\leq i,j\leq p} ,\mbox{ where  } m_{ii} =
  H_{\varepsilon}(x_i ,x_i) , m_{ij} = -G_{\varepsilon}(x_i ,x_j) ,i\neq j
\end{eqnarray}
\noindent
and define $\rho_{\varepsilon}(x) $ as the least eigenvalue of $M$  ($ \rho_{\varepsilon}(x) = -\infty $ if $x_i = x_j$ for some  $i\neq j$). \\
For $a\in\R^3$ and $\lambda >0,$  $\delta_{(a,\lambda )}$ denotes
the function
\begin{eqnarray}\label{e:24}
 \delta _ {(a,\lambda )}(x) = c_0 \frac{\lambda ^{1/2}}{(1+
 \lambda^2|x-a|^2)^{1/2}}.
\end{eqnarray}
\noindent It is well known (see \cite{CGS}) that if $c_0$ is suitably chosen $ (c_0
= 3^{1/4})$, $\delta _{(a,\lambda )}$ are the only solutions of
\begin{eqnarray}\label{e:25}
\quad -\Delta u =  u^{5} , u>0 \mbox{  in  } \R^3
\end{eqnarray}
\noindent
and they are also the only minimizers for the Sobolev inequality
\begin{eqnarray}\label{e:26}
 S =\inf\{|\nabla u|^{2}_{L^2(\R^3)}|u|^{-2}_{L^{6}(\R^3)}
,\, s.t.\, \nabla u\in L^2 ,u\in L^{6} ,u\neq 0 \}.
\end{eqnarray}
\noindent We also denote by $P_\e\delta_{(a,\lambda)}$ the
projection of $\delta_{(a,\lambda )}$ on $H_0^1(A_{\varepsilon}),$
that is,
$$
\Delta P_\e\delta_{(a,\lambda )} = \Delta \delta_{(a,\lambda )}
\mbox{  in  } A_{\varepsilon},  P_\e\delta_{(a,\lambda )}  = 0
\mbox{ on } \partial A_{\varepsilon},
$$
and by $\th_{(a,\l)}= \d_{(a,\l)}- P_\e\d_{(a,\l)}$. We define on
$H_0^1(A_{\varepsilon})\setminus \left\{0\right\}$ the functional
\begin{eqnarray}\label{e:17}
\quad J_{\varepsilon}(u) = \frac{\int_{A_{\varepsilon}}|\nabla u|^2}{\left(
 \int_{A_{\varepsilon}} u ^6 \right)^{1/3}}
\end{eqnarray}
\noindent whose positive critical points, up a multiplicative
constant, are solutions  of $(P_\e)$. Lastly, let
\begin{eqnarray*}
\langle u,v\rangle = \int_{A_\e}\n u\n v, \quad ||u||=
\left(\int_{A_\e}|\n u|^2 \right)^{1/2},\quad u,\,v\, \in
H^1_0\left(A_\e\right).
\end{eqnarray*}
\section{Asymptotic behavior of bounded energy solutions}
This section is devoted to the study of  the asymptotic behavior
of bounded energy solutions of $(P_\e )$. Such a precise
description is  cornerstone in the proof of our results. It says,
roughly speaking, that our solutions concentrate  at a finite
number of points such that the
distance of one of them to the other is at least comparable  to
$\e$.

In  the sequel of this paper we consider a solution $u_\e$ of
$(P_\e)$ which satisfies \be\label{borne} \int_{A_\e} \mid \n
u_\e\mid^2 \leq C, \ee where $C$ is a positive constant
independent of $\e$. Our aim in this section is to prove the
following result:
\begin{thm} \label{t:21}
Let $u_{\varepsilon}$ be a solution of problem $(P_\varepsilon)$
which satisfies \eqref{borne}. Then, after passing to a
subsequence, there exist  $p\in\N^*, (x_{1,
\varepsilon},...,x_{p,\varepsilon})\in A_\varepsilon^p ,
(\lambda_{1, \varepsilon},...,\lambda_{p, \varepsilon})\in
(\R_+^*)^p$, and a positive constant $\a > 0$ such that:
\begin{align*}
&\biggl|\biggl|u_\varepsilon - \sum_{i=1}^{p} P_{\varepsilon}
\delta_{(x_{i , \varepsilon},\lambda_{i
,\varepsilon})}\biggr|\biggr| \rightarrow 0,\,
\lambda_{i,\varepsilon}d_{i,\varepsilon}\rightarrow +\infty \,{
for }\, 1\leq i\leq p\quad
\mbox{as}\quad \e\to 0,\\
&\l_{i,\e}\mid x_{i,\e} - x_{j,\e}\mid \to \infty  \mbox{ as }
\e\to 0, \mid x_{i,\e} - x_{j,\e}\mid  \geq \a \, \e
\,\,\mbox{for}\,\, i\ne j,
\end{align*}
where $ d_{i,\varepsilon}=d(x_{i,\varepsilon},\partial A_\e)$ and
$\l_{i,\e}= 3^{-1/2}\left(u_\e(x_{i,\e})\right)^{2}$.
\end{thm}

\begin{rem}
The above Theorem is true in all dimensions $n \geq 3$, however a
weaker version used in \cite{BEH} was enough to derive the
equivalent of our result in dimension $n \geq 4$.
\end{rem}

To prove Theorem \ref{t:21}, we  start by establishing some
useful facts. Let $x_{1,\e} \in A_\e$ be such that
$$
u_\e(x_{1,\e}) = \max_{A_\e} u_\e := M_{1,\e}.
$$
Let $\widetilde A_\varepsilon =
M_{1,\varepsilon}^{2}(A_\varepsilon - x_{1,\varepsilon })$, and
denote by $v_\varepsilon$ the function defined on $\widetilde
A_\varepsilon$ by
\begin{eqnarray}
v_\varepsilon(y) = M_{1,\varepsilon}^{-1}u_\varepsilon (x_{1,
\varepsilon} + M_{1,\varepsilon}^{-2} y).
\end{eqnarray}
By Lemma 2.3 of \cite{BEH}, we know that:
$$
M_{1,\e}^{2}d(x_{1,\e}, \partial A_\e) \to +\infty
\quad\mbox{as}\quad \e\to 0.
$$
Furthermore, $v_\e \to \d_{(0,\a_0)}$ in $C^2_{loc}(\R^3)$ as
$\e\to 0$, where $\a_0=3^{-1/2}$.\\ Now, we prove the following
crucial lemma:
\begin{lem}\label{l:1}
There exist positive constants $\d$ and ${\bar c}$ such that
$$
\max_{\mid y\mid \leq \d\e M_{1,\e}^{2}} \mid
v_\e(y)-\d_{(0,\a_0)}(y)\mid \leq {\bar c}\left(\e
M_{1,\e}^2\right)^{-1}.
$$
\end{lem}
\begin{pf}
First, it follows from Lemma 3.2 of \cite{CL}, that there exist
positive constants $\d$ and ${\bar c}$ such that \be\label{(*)}
v_\e(y) \leq {\bar c} \d_{(0,\a_0)}(y) \quad\mbox{for}\quad \mid
y\mid \leq \d \e M_{1,\e}^{2}. \ee
 Now, let
$$
m_\e= \max_{\mid y\mid \leq \d\e M_{1,\e}^{2}} \mid
v_\e(y)-\d_{(0,\a_0)}(y)\mid:= \mid
v_\e(y_\e)-\d_{(0,\a_0)}(y_\e)\mid.
$$
Arguing by contradiction, we assume that
$m_\e \e M_{1,\e}^2 \to +\infty$ as $\e\to 0$.\\
Let $w_\e(y)= m_\e^{-1} \left(v_\e(y)-\d_{(0,\a_0)}(y)\right)$,
$w_\e$ satisfies
$$
\D w_\e + f_\e w_\e = 0\quad\mbox{with}\quad f_\e = \frac{v_\e^{5}
- \d_{(0,\a_0)}^{5}}{v_\e - \d_{(0,\a_0)}}.
$$
 By
\eqref{(*)}, we have \be\label{(**)} \mid f_\e\mid \leq c\left(1+
\mid y\mid\right)^{-4}\quad \mbox{for}\quad \mid y\mid \leq \d \e
M_{1,\e}^{2}. \ee Applying the Green's representation leads to
$$
w_\e(y)= a\left(\int_{B_\e} G_{B_\e}(y,x)f_\e(x)w_\e(x)dx -
\int_{\partial B_\e} \frac{\partial G_{B_\e}}{\partial
\nu}(y,x)w_\e(x)d\sigma(x)\right),
$$
where $a=\left(meas\left(S^2\right)\right)^{-1}$, $B_\e= B(0, \d\e
M_{1,\e}^{2})$, $\nu$ is the outward normal to $\partial B_\e$ and
$G_{B_\e}$ is the Green's function of $\D$ under Dirichlet
boundary conditions in $B_\e$. Using \eqref{(*)} and \eqref{(**)}
yields
\begin{align}\label{(***)}
\mid w_\e(y)\mid &\leq c \int_{B_\e}\frac{dx}{\mid y-x\mid
\left(1+\mid x\mid\right)^4} + \frac{c}{m_\e \d \e M_{1,\e}^{2}}\notag\\
&\leq c\left( \left(1+ \mid y\mid\right)^{-2} + \left( m_\e \d\e
M_{1,\e}^{2}\right)^{-1}\right).
\end{align}
It follows that $w_\e$ is bounded and by elliptic standard
estimates $w_\e$ converges, up to some subsequence, in the
$C^2_{loc}$-norm to a function $w$ satisfying
 \be\label{E}
\begin{cases}
\D w + 5\d_{(0,\a_0)}^{4}(y)w(y)=0\quad\mbox{in}\quad \R^3\\
\mid w(y)\mid \leq c\left(1+\mid y\mid \right)^{-2}.
\end{cases}
\ee
By Lemma 2.4 of \cite{CL}, every solution of \eqref{E} can be written as
$$
w(y)=\sum_{j=1}^3 a_j \frac{\partial\d_{(0,\a_0)}}{\partial y_j} +
a_0 \left( y \cdot \, \n \d_{(0,\a_0)}(y)+\frac{1}{2}
\d_{(0,\a_0)}(y)\right)
$$
for some constants $a_j\geq 0$, $j=0,...,3$. Since
$w(0)=\frac{\partial w}{\partial y_j}(0)=0$, we obtain that
$a_j=0$ for $0\leq j\leq 3$, namely, $w\equiv 0$. Since
$w_\e(y_\e)=1$, it follows that $\mid y_\e\mid \to +\infty$ as
$\e\to 0$. Applying \eqref{(***)} at $y=y_\e$ gives \be\label{4*}
1=\mid w_\e(y_\e)\mid \leq c \left( \left(1+\mid y_\e\mid
\right)^{-2} +\left( m_\e \d \e M_{1,\e}^{2}\right)^{-1}\right).
\ee Since the right hand-side of \eqref{4*} goes to zero, as
$\e\to 0$, we derive a contradiction. Thus $m_\e \e \d M_{1,\e}^2$
must be bounded and the proof of our lemma follows.
\end{pf}
\begin{lem}\label{l:2}
Let $\d$ be the positive constant stated in Lemma \ref{l:1}. Then we have
$$
\int_{B(x_{1,\e},\d\e)} u_\e^{6}= S_3 + o(1) \quad\mbox{as}\quad
\e\to 0,
$$
where $S_3=S^{3/2}$ and $S$ is the Sobolev constant defined in
\eqref{e:26}.
\end{lem}
\begin{pf}
We have
\begin{align*}
\int_{B(x_{1,\e},\d\e)} u_\e^{6}&=\int_{B(0,\d\e M_{1,\e}^{2})} v_\e^{6}\\
&=  \int_{B(0,\d\e M_{1,\e}^{2})} \d_{(0,\a_0)}^{6} + O\left(
\int_{B(0,\d\e M_{1,\e}^{2})}
\d_{(0,\a_0)}^{5}|v_\e-\d_{(0,\a_0)}|+
|v_\e-\d_{(0,\a_0)}|^6\right)\\
&=  \int_{B(0,\d\e M_{1,\e}^{2})} \d_{(0,\a_0)}^{6} + O\left( \mid
v_\e-\d_{(0,\a_0)}\mid_{L^{6}(B(0,\d\e M_{1,\e}^{2}))}\right).
\end{align*}
Using Lemma \ref{l:1} and the fact that $\e M_{1,\e}^{2} \to
+\infty$ as $\e\to 0$, we easily derive our lemma.
\end{pf}

Now, we are in the position to prove Theorem \ref{t:21}.\\
\begin{pfn}{Theorem \ref{t:21}}
We distinguish two cases:\\
{\bf Case 1.}.  $\int_{A_\e}\mid u_\e - P_\e\d_{(x_{1,\e},
\l_{1,\e})}\mid^{6} \to 0$ as $\e \to 0$, where
$\l_{1,\e}=\a_0M_{1,\e}^{2}$. In this case we are done, the number
of blow up points in the Theorem is reduced to $1$,
that is, $p=1$.\\
{\bf Case 2.}.  $\int_{A_\e}\mid u_\e - P_\e\d_{(x_{1,\e},
\l_{1,\e})}\mid^{6} \not\to 0$ as $\e \to 0$.  We are going to
study this case. First, let us prove that \be\label{second1}
\int_{A_\e\setminus B(x_{1,\e}, \d\e)} u_\e^{6} \not\to
0\quad\mbox{as}\quad \e\to 0. \ee Observe that
 \be\label{second2}
\int_{A_\e\setminus B(x_{1,\e}, \d\e)} P_\e\d_{(x_{1,\e},
\l_{1,\e})}^{6} \leq \int_{A_\e\setminus B(x_{1,\e}, \d\e)}
\d_{(x_{1,\e}, \l_{1,\e})}^{6} = \int_{\tilde{A}_\e\setminus B(0,
\d\e M_{1,\e}^{2})} \d_{(0, \a_0)}^{6}\to  0\mbox{ as } \e\to 0.
  \ee
where we have used the fact that $\e M_{1,\e}^{2} \to \infty$ and
$\d_{(0,\a_0)} \in L^{6}(\R^3)$.\\ By Lemma \ref{l:1} and the fact
that $\e M_{1,\e}^2\to \infty$, it is easy to derive
 \be\label{e:}
 \int_{B_\e} |u_\e-P\d_{(x_{1,\e},\l_{1,\e})}|^6 \to 0\quad \mbox{as}\quad \e \to
 0.
 \ee
Clearly, \eqref{second2} and \eqref{e:} imply \eqref{second1}.
Now, we set
$$
u_\e(x_{2,\e})=\max_{A_\e\setminus B(x_{1,\e},\d\e)} u_\e:= M_{2,\e}.
$$
It is clear that $|x_{1,\e}-x_{2,\e}|\geq \d \e$.\\ By
\eqref{second1}, there exists $c>0$ such that
$$
c \leq \int_{A_\e\setminus B(x_{1,\e}, \d\e)} u_\e^{6} \leq
M_{2,\e}^{4} \int_{A_\e}u_\e^2(x)dx.
$$
But, we have
$$
 \int_{A_\e}u_\e^2(x)dx = \e^3\int_{D_\e} \tilde{u}_\e ^2(X)dX\leq
 \frac{\e^3}{c_\e}\int_{D_\e}\mid\n\tilde{u}_\e(X)\mid^2dX=
 \frac{\e^2}{c_\e}\int_{A_\e}\mid\n u_\e(x)\mid^2dx\leq \frac{C\e^2}{c_\e},
$$
where $\tilde{u}_\e(X)=u_\e(\e X)$, $D_\e=\var(A_\e)$, with $\var:
x\mapsto \var(x)=\e^{-1}x$, and $c_\e>0$. By Lin \cite{Lin}, we
have $c_\e \to c>0$ as $\e\to 0$. We derive that $\e M_{2,\e}^{2}
\not\to 0$ as $\e\to 0$ and therefore as in Lemma 2.3 of
\cite{BEH}, we have that $ M_{2,\e}^{2} d(x_{2,\e}, \partial
A_\e)\to +\infty$ as $\e\to 0$. This implies that $
M_{2,\e}^{2}\mid x_{1,\e}-x_{2,\e}\mid \to +\infty$ as $\e\to 0$.
Now, for $y\in E_\e:= M_{2,\e}^2\left(A_\e - x_{2,\e}\right)$, we
introduce the following function
$$
U_\e(y)= M_{2,\e}^{-1}u_\e\left(x_{2,\e}+ M_{2,\e}^{-2}y\right).
$$
It is easy to see that $U_\e$  is bounded by 1 in
$B(0,(1/2)M_{2,\e}^2\mid x_{2,\e}-x_{1,\e}\mid)$. Therefore, $U_\e
\to \d_{(0,\a_0)}$ in $C^2_{loc}(\R^3)$ as $\e\to 0$. Thus, we
have obtained in Case 2 a second blow up point. It is clear that
we can iterate such a process. But, since the energy of $u_\e$ is
bounded such a process stops after finitely steps, and the proof
of our Theorem is thereby completed.
\end{pfn}

\section{Proof of Theorem \ref{t:11}}
This section is devoted to the proof of Theorem \ref{t:11}. To
this aim, we first study the location of blow up points that we
found in Section 3. To this goal, we need a rather delicate
analysis and careful estimates. First, we start by the general
setting. Let, for $p\in\N^*$ and $\eta >0$ given
\begin{align*}
V_\varepsilon (p,\eta)=\biggl\{& u\in \Sigma ^{+}(A_\varepsilon
)\, \mbox{ s.t } \exists \, y_1,..., y_p \in A_\varepsilon,
\exists \, \lambda_1,...,\lambda_p  >\frac{1}{\eta}
\, \, \mbox{with} \\
  & \biggl|\biggl|u -C(p) \sum_{i=1}^pP_\varepsilon\delta_{( y_i,\lambda_i)}
\biggr|\biggr|  <\eta ,\,  \lambda_i \, d(y_i ,
\partial A_\varepsilon )> \frac{1}{\eta }\, \forall \, i,
  \quad \varepsilon_{ij} <\eta\, \forall i\neq j  \biggr\},
\end{align*}
\noindent where $\Sigma ^+ (A_{\varepsilon }) = \{ u \in H_{0}^{1}
 (A_\varepsilon) /\, u > 0,\, \, ||u|| = 1\}$ and
 $\e_{ij}=(\l_i/\l_j+\l_j/\l_i+\l_i\l_j|y_i-y_j|^2)^{-1/2}$.\\
If a function $u$ belongs to $V_\varepsilon(p,\eta ),$ then, for
$\eta >0$ small enough, the minimization problem
\begin{eqnarray}\label{e:212}
\min_{
 \alpha_i ,\lambda_i >0,\  y_i\in A_\varepsilon
} \biggl|\biggl|u - \sum_{i=1}^p\alpha_iP_\varepsilon\delta_{(
 y_i,\lambda_i)}\biggr|\biggr|
\end{eqnarray}
\noindent
has a unique solution, up to permutation (see Lemma  A.2 in \cite {BC}).\\
Therefore, for $\varepsilon >0$ sufficiently small, $u_\e $
(solution of $(P_\e)$) can be uniquely written as
\begin{eqnarray}\label{e:213}
 \widetilde{u}_\varepsilon := \frac{u_\varepsilon}{||
u_\varepsilon||} = \sum_{i=1}^p \alpha_{i,\varepsilon}
P_\varepsilon\delta_{(x_{i,\varepsilon} , \lambda_{i,\varepsilon}
)} + v_\varepsilon,
\end{eqnarray}
\noindent where $v_\varepsilon$ satisfies the following conditions:
$$
(V_0) \, \,  \langle v_\varepsilon , P_\varepsilon
\delta_{(x_{i,\varepsilon} ,\lambda_{i,\varepsilon} )} \rangle =
\langle  v_\varepsilon ,\frac{\partial P_\varepsilon
\delta_{(x_{i,\varepsilon},\lambda_{i,\varepsilon} )}}{\partial
\lambda_{i,\e}}\rangle =\langle v_\varepsilon ,\frac{\partial
P_\varepsilon \delta_{(x_{i,\varepsilon}
 ,\lambda_{i,\varepsilon} )}}{\partial (x_{i,\e})_k}\rangle =
 0\quad \forall \, i,
$$
where $(x_{i,\e})_k$ is the $k$th component of $x_{i,\e}$, $k\in
\{1,2,3\}$ and $\a_{i,\e}$ satisfies :
$$
 J(u_{ \e})^3 \alpha_{j,\e}^{4} = 1 +o(1)\, \,
\forall j.
$$
\noindent To simplify the notations,  we write $\alpha_i,x_i,
 \lambda_i, \delta_i $,  $P\delta_i$  and   $\th_i$ instead of $\alpha_{i,\varepsilon} ,
x_{i,\varepsilon} , \lambda_{i,\varepsilon},
\delta_{(x_{i,\varepsilon} ,\lambda_{i,\varepsilon})} $,
$P_\e\delta_{(x_{i, \varepsilon} ,\lambda_{i,\varepsilon})}$ and
$\th_{(x_{i, \varepsilon} ,\lambda_{i,\varepsilon})}$ respectively
and we also write $u_\varepsilon$ instead of
$\widetilde{u}_\varepsilon.$

As a consequence of Theorem \ref{t:21}, it is easy to obtain the
following result
\begin{cor}\label{c:A} For each $i$, we denote by $B_i:=B(x_i,\a
d_i /4)$. For $i\ne j$, we have
$$(a) \quad \e_{ij}\leq \frac{c}{(\l_id_i\l_jd_j)^{1/2}},\quad (b)
\quad \l_i\frac{\partial\e_{ij}}{\partial\l_i}=
-\frac{1}{2}\e_{ij}(1+o(1)),\quad (c)\quad B_i\cap
 B_j=\emptyset.$$
\end{cor}
\begin{pf}
The proof is immediate since $|x_i-x_j|\geq \a \e$ for each $i\ne
j$ and $d_i \leq \e$ for each $i$.
\end{pf}

 Now, let us recall the estimate of the $v_\varepsilon$-part
of $u_\e$.
\begin{pro}\cite{BEH}\label{p:24}
Let $v_\varepsilon$ be defined by \eqref{e:213}. Then, we have the
following estimate
$$
||v_{\varepsilon}|| \leq c \sum_i\frac{1}{\lambda_id_i}   +
c\sum_{i\ne j}\varepsilon_{ij}\left( Log \varepsilon_{ij}^{-1}
\right)^{1/3}.
$$
\end{pro}

In the next propositions, we give  useful expansions of the
gradient of $J$ which allows us to characterize the concentration
points given by Theorem \ref{t:21}.

Regarding the estimate of $||v_\e||^2$, it is negligible with
respect to the principle part of Proposition 3.2 of \cite{BEH},
however it is of the same order as the principle part of
Proposition 3.3 of \cite{BEH}. Following an idea introduced by O.
Rey \cite{R2} and the fact that the balls $B_i$ are disjoints, we
are able to improve the terms which contain $v_\e$ and
therefore we can obtain the analogue of Proposition 3.3 of
\cite{BEH}.

\begin{pro}\label{p:214}
 For each $i$, we have the following expansion
\begin{align*}
\langle \nabla  J(u_\varepsilon ),\lambda_i{\frac{
\partial P\delta_i}{\partial\lambda_i}}\rangle
 &=  2J(u_\varepsilon )c_1\biggl(-\frac{\a_i}{2}
\frac{H_\varepsilon (x_i,x_i)}{\lambda_i }(1+o(1))\\
 &- \sum_{j\ne i}\alpha _j \left(\lambda _i \frac{\partial
\varepsilon_{ij}}{\partial\lambda_i}  + \frac{1}{2}
\frac{H_\varepsilon (x_i,x_j)}{(\lambda_i\lambda_j)^{1/2} }
\right)(1+o(1)) +R \biggr),
\end{align*}
\noindent where $c_1$ is a positive constant and $ R =
O\left(\sum_1^p{(\lambda_kd_k )^{-2}} + \sum_{k\neq
r}^{}\varepsilon_{kr}^2\left(Log
\varepsilon_{kr}^{-1}\right)^{2/3}\right)$.
\end{pro}
\begin{pf}
It follows from Lemma \ref{l:A0}, Proposition \ref{p:24} and the
fact that $v_\e$ satisfies $(V_0)$.
\end{pf}

\begin{pro}\label{p:215}
For each $i$, we have the following expansion
\begin{align*}
\langle\nabla J(u_\varepsilon ),\frac{1}{\lambda_i}\frac{
\partial P\delta_i}{\partial x_i}  \rangle = &
J(u_\e )c_1 \biggl( -2 \sum_{j\ne i}\alpha _j \left(
\frac{1}{\lambda_i} \frac{\partial \varepsilon_{ij}}{\partial x_i}
 - \frac{1}{\lambda_i(\lambda_i \lambda_j)^{1/2}}
\frac{ \partial H_\varepsilon(x_i ,x_j )}{\partial x_i} \right)\\
 & +\frac{\alpha_i}{\lambda_i^{2}}\frac{ \partial
H_\varepsilon(x_i ,x_i )}{\partial x_i} + o\biggl(\sum_1^p
\frac{1}{\left(\l_kd_k\right)^2}\biggr) \biggr).
\end{align*}
\noindent
\end{pro}
\begin{pf}
It follows from Lemmas \ref{l:A1}, \ref{l:A2}, \ref{l:A5},
\ref{l:210}, \ref{l:213}, Proposition \ref{p:24} and the fact that
$v_\e$ satisfies $(V_0)$. The negligible terms which appear in
those estimates can be written as $o((\l_1d_1)^{-2})$ since
$|x_i-x_j|\geq \a \e$ for each $i\ne j$ and $d_k \leq \e$ for each
$k$.
\end{pf}

 Now, we order all the $\lambda_id_i$'s : $\lambda_{1}d_{1}\leq
\lambda_{2}d_{2} \leq ... \leq\lambda_{p}d_{p}$.

First, we introduce the set of  indices $i$ such that $\l_id_i$
and $\l_1d_1$ are of the same order. Let $C_1$ be a large positive
constant and define
 \be \label{l}
 I=\{1\}\cup \{ i / \l_k d_k \leq C_1
\l_{k-1}d_{k-1} \mbox{ for each } k\leq i\} := \{1,2,...,l\}. \ee

 Secondly, we define a subset of $I$ such that the distance between
the points is at most comparable to their distances to the
boundary. Let $C_0$ be a large positive constant, we define
 \be \label{e:323}
 B= \left\{i\in I / \exists k_1,...,k_m \in I \,
s.t.\, \,  k_1 = i,...,k_m = 1\,  ;\,  {|x_{k_j}
-x_{k_{j+1}}|}\leq C_0\min(d_{k_j} ,d_{k_{j+1}})\right\}.
 \ee

\begin{lem}\label{l:33}
Let B be defined by \eqref{e:323}. Then, $\{1\}\varsubsetneq B$.
\end{lem}
\begin{pf}
First, we remark that Proposition \ref{p:214} implies immediately
that $p\geq 2$. To prove our lemma, we argue by contradiction.
We assume that $ B =\{ 1\}.$  \\
Using Proposition \ref{p:214}, and the fact that
$H_\varepsilon(x_i ,x_i) \sim{c}/{d_i}$ (see \cite{AE1}), we derive
\begin{eqnarray}\label{e:324}
0 = \langle\n J(u_\e ),\l _1\frac{\partial P\d _1}{\partial\l
_1}\rangle \leq -\frac{c}{(\lambda_1d_1)} + O( \sum_{k\neq
1}\varepsilon_{1k}).
\end{eqnarray}
\noindent Two cases may occur. If $k>l$ where $l$ is defined by
\eqref{l}, then by Corollary \ref{c:A}, we have
$$
 \varepsilon_{1k} \leq
\frac{c}{(\lambda_kd_k\lambda_1d_1)^{1/2}} \leq\frac{1}{C_1^{1/2}}
\frac{1}{((\lambda_ld_l)(\lambda_1d_1))^{1/2}}=
o\left(\frac{1}{\l_1d_1}\right) \quad ( \mbox{ for } C_1 \mbox{
large enough}).
$$
In the other case, we have $|x_1 -x_k|\geq C_0 \mbox{ min }
(d_1,d_k)$, then
$$
\e_{1k} \leq \left(\frac{1}{\lambda_1 \lambda_k|x_1 - x_k|^2}
\right)^{1/2} \leq \frac{2}{C_0^{1/2}}\frac{1}{\left((
\lambda_1d_1)(\l_kd_k) \right)^{1/2}}
 = o\left(\frac{1}{\l_1d_1}   \right)\quad ( \mbox{ for } C_0 \mbox{
large enough}).
$$
Thus \eqref{e:324} yields a contradiction and the result follows.
\end{pf}

 Next, our goal is to prove the following crucial result:
\begin{pro}\label{p:p1}
Let $x_{1,\e}$,..., $x_{p,\e}$ be the points given by Theorem
\ref{t:21}. Then, we have $p\geq 2$ and there exist
$k\in\{2,..,p\}$, $i_1$,..., $i_k \in\{1,...,p\}$ such that
$$
d\rho_\e(x_{i_1,\e},...,x_{i_k,\e}) \to 0 \quad \mbox{ and }\quad
d^2 \n \rho_\e(x_{i_1,\e},...,x_{i_k,\e}) \to 0,\quad  \mbox{ as
}\e\to 0,
$$
where $d=\min_{1\leq r\leq k}d(x_{i_r,\e}, \partial A_\e)$. In
addition, we have $\forall m,\,r\,\in\{1,...,k\}$
$|x_{i_m,\e}-x_{i_r,\e}| \leq C_0' d$, where $C_0'$ is a positive
constant independent of $\e$.
\end{pro}
\begin{pf}
Let $k=\mbox{card }B$ that is $B=\{i_1,...,i_k\}$. By Lemma
\ref{l:33}, we have $k\geq 2$.\\ Let $M_B=(m_{ij})_{i,j\in B}$ be
the matrix defined by \eqref{e:23} and let
$\rho_B=\rho_\e(x_{i_1,\e},...,x_{i_k,\e})$ be the least
eigenvalue associated to $M_B$. We denote by $e$ the eigenvector
associated to $\rho_B$ whose norm is 1. We know that all
components of $e$ are strictly positive (see \cite{BLR}). Let
$\eta >0$ be such that for any $\gamma$ belongs to a neighborhood
$C(e,\eta )\subset \{ y \in (\R^*_+)^k s.t.\, \,  ||y|^{-1} y - e|
< \eta \}$, we have
\begin{eqnarray}\label{e:326}
   ^T\gamma M_B\gamma  - \rho_B|\gamma |^2 \leq
 \frac{c_2}{d}|\gamma |^2 \mbox{ and } ^T
\gamma\frac{\partial M_B}{\partial x_i}\gamma =
 \left(\frac{\partial \rho _B}{\partial x_i} +
 o(\frac{1}{d^2})  \right)|\gamma |^2
 \end{eqnarray}
and for $\gamma\in (\R^*_+)^k\setminus C(e, \eta )$ , we have
\begin{eqnarray}\label{e:327}
 ^T\gamma M_B\gamma - \rho_B|\gamma |^2\geq
{c_3|\gamma |^2}d^{-1}.
\end{eqnarray}
First, we study the vector $\Lambda$ defined by $ \Lambda  =
\left({\lambda_{i_1}^{-1/2}},
...,{\lambda_{i_k}^{-1/2}}\right)$.\\
{\bf Claim 1.} We have
$\Lambda \in C(e, \eta ).$\\
\noindent {\bf Proof of Claim 1.} We argue by contradiction.
Assume that $\Lambda\in (\R^*_+)^k\setminus C(e,\eta )$. Let
$$
\Lambda  (t) = |\Lambda |\frac{(1-t)\Lambda  +
 t|\Lambda  |e}{| (1-t)\Lambda + t|\Lambda |e |}
: =\frac{y(t)}{|y(t)|}.
$$
\noindent From Proposition \ref{p:214}, we derive
$$
\langle\n J(u_\e ),Z\rangle_{|t=0} = -c\frac{d}{dt}
 \left(^T\Lambda  (t)M_B\Lambda (t)  \right)  +
 O\left( \sum_{i\in B,j\notin  B}
\varepsilon_{ij} \right)+o\left(  \frac{1} {\lambda_1d_1} \right)
$$
where $Z$ is the vector field defined on the variables $\l $ along
the flow line defined by $\L (t)$.\\
Observe that
\begin{align*}
\frac{d}{dt} & \left( ^T\Lambda  (t)M_B\Lambda  (t) \right)
=\frac{d}{dt}\left(\frac{ ^T\Lambda  (t)M_B \Lambda (t)}{|\Lambda
(t) |^2} |\Lambda (0)|^2
 \right)\nonumber\\
&= |\Lambda (0)|^2\frac{d}{dt}\left( \rho_B
+\frac{(1-t)^2}{|y(t)|^2}(^T\Lambda  (0)M_B\Lambda (0)
 - \rho _B|\Lambda (0) |^2  ) \right)\nonumber\\
&= |\Lambda (0)|^2\left( \frac{2(1-t)}{|y(t)|^4}(^T \Lambda
(0)M_B\Lambda (0)  - \rho_B|\Lambda (0) |^2
 )(-(1-t)|\Lambda (0)|<e, \Lambda (0)> -t|\Lambda |^2) \right).
\end{align*}
\noindent Thus
\begin{align*}
\langle\n J(u_\e ),Z\rangle_{|t=0} =& -\frac{2c}{|\Lambda |^2}
(^T\Lambda M_B\Lambda   - \rho _B|\Lambda  |^2  )
 (-|\Lambda |<e, \Lambda (0)>) \\
 &+ o\left(\frac{1}{
\lambda_1d_1 }  \right)  + O\biggl( \sum_{i\in B, j\notin
 B}\varepsilon _{ij} \biggr).
\end{align*}
\noindent Since $|e| = 1,$  then there exists $m$ such that
$e_{i_m} \geq \frac{1}{k}.$  Thus
$$
<e,\Lambda (0)> = \sum_j e_{i_j} \Lambda _{i_j} \geq
\frac{1}{k}\Lambda _{i_m}.
$$
\noindent Using \eqref{e:327}, we obtain
\begin{eqnarray*}
\langle\n J(u_\e ),Z\rangle_{|t=0}&\geq & \frac{cc_3}{d}| \Lambda
|\Lambda _{i_m} + o\left(\frac{1}{\lambda_1d_1}
 \right)  + O\biggl( \sum_{i\in B,j\notin  B}
\varepsilon _{ij} \biggr)\nonumber \\
&\geq & \frac{c}{\left( \lambda _1 d_1 \lambda _{i_m} d_{i_m}
 \right)^{1/2}} + o\left ( \frac {1}{
\lambda _1 d_1 }\right) + O \biggl( \sum _{i \in B, \, j \notin
 B} \varepsilon _{ij} \biggr).
\end{eqnarray*}
\noindent As in the proof of Lemma \ref{l:33}, we have
 \be \label{e:C1}
\varepsilon_{ij} = o\left( \frac{1}{(\lambda_1
d_1\lambda_{i_m}d_{i_m} )^{1/2}} \right)\quad \forall \, i\in B,
\, \forall \, j\notin B.
 \ee
 Thus
$$
 0\geq  \frac{c}{(\lambda_1d_1\lambda_{i_m}d_{i_m}
)^{1/2}}  +o\left( \frac{1}{ \lambda_1d_1 }  \right) \geq
\frac{1}{\lambda_1d_1 }  \left(
 \frac{c}{C_1^{k/2}} +o(1)   \right) >0.
$$
This yields a contradiction and our claim follows.

Now, we will prove that
 \be\label{e:AA}
 d\rho _B \longrightarrow 0 , \mbox{ as
}\varepsilon \longrightarrow 0.
 \ee
Using Proposition \ref{p:214} and \eqref{e:C1}, we have
\begin{align}
0 & =  \sum_{i\in B}
\langle\n J(u_\e ),\l _i\frac{\partial P\d _i}{\partial \l _i}\rangle \notag\\
 & = \sum_{i\in B}\left[\frac{
H_\varepsilon (x_i ,x_i )}{\lambda_i}(1+ o(1)) - \sum_{j\ne i,
j\in B}(\varepsilon_{ij} -  \frac{ H_\varepsilon (x_i ,x_j
)}{(\lambda_i \lambda_j)^\frac{1}{2}})(1+ o(1)) +O(\sum_{j\notin
B} \varepsilon_{ij} )+R \right]\notag\\
 & =  ^T \Lambda M_B\Lambda +o\left( \frac{1}
{\lambda_1d_1 } \right).\label{e:328}
\end{align}
We assume, arguing by contradiction, that $d\rho _B \not
\longrightarrow 0$, when $\varepsilon \longrightarrow 0.$
Therefore, there exists $C_4>0$ such that $d|\rho _B|\geq C_4$.
Now, we distinguish two cases \\
{\bf $1^{st}$case:} $\rho _B>0$. In this case, we derive from
\eqref{e:328}
$$
0\geq \rho _B|\Lambda |^2 + o\left(\frac{1}{ \lambda_1d_1 }
\right) \geq C_4\frac{ |\Lambda |^2}{d} +o\left(\frac{1}{
\lambda_1d_1 } \right) >0.
$$
\noindent This yields a contradiction. \\
{\bf {$ 2^{nd}$ case:}}  $\rho_B < 0. $ In this case, using Claim
1, we derive from \eqref{e:326} and \eqref{e:328},
\begin{eqnarray*}
0\leq \rho _B|\Lambda |^2 +\frac{c_2 | \Lambda |^2}{d}
+o\left(\frac{1}{ \lambda_1d_1 } \right)& \leq&
 \frac{|\Lambda |^2}{d}\left(
 \rho_B \, d +c_2   \right) +o\left(
\frac{1}{\lambda_1d_1 } \right)\nonumber \\
& \leq& \frac{|\Lambda |^2}{d}(-C_4 +c_2)
+o\left(\frac{1}{\lambda_1d_1 } \right).
\end{eqnarray*}
\noindent If we choose $c_2 \leq \frac{1}{2}C_4,$ we obtain a
contradiction. Thus, \eqref{e:AA} follows.

In order to complete the proof of Proposition \ref{p:p1},  it
remains to prove that:
 \be \label{e:AB}
 d^2 \n \rho_B \to 0, \quad \mbox{ as } \e \to 0.
 \ee

We assume, arguing by contradiction, that $d^2\nabla \rho_B \not
\longrightarrow 0$  when  $\varepsilon\longrightarrow 0.$\\
For $i\in B$, using Proposition \ref{p:215}, we derive
\begin{eqnarray*}
0  = ^T\Lambda \frac{\partial M_B}{\partial x_i} \Lambda +
O\left(\sum_{j\notin B}\frac{
\partial\varepsilon_{ij}}{\partial x_i} -
 \frac{1}{( \lambda_i\lambda_j )^\frac{1}{2}}
\frac{\partial H_\varepsilon}{\partial x_i}(x_i ,x_j )\right)
+o\left( \frac{1}{d_i}\frac{1} {(\lambda_1d_1 )}\right).
\end{eqnarray*}
 Observe that $|\partial H/\partial x_i (x_i,x_j)| \leq c (d_i
 |x_i-x_j|)^{-1}$. Thus, as in the proof of Lemma \ref{l:33}, we
 prove that, for $i\in B$ and $j\notin B$,
\begin{eqnarray*}
\left|\frac{\partial\varepsilon_{ij}}{\partial x_i}\right| +
 \frac{1}{( \lambda_i\lambda_j )^{1/2}}\left|\frac{
\partial H_\varepsilon}{\partial x_i}(x_i,x_j)\right|
= o\left(\frac{1}{d(\lambda_1d_1)}  \right).
\end{eqnarray*}
Therefore, by \eqref{e:326}, we have
$$
0= ^T\Lambda \frac{\partial M_B}{\partial x_i}\Lambda +
o\left(\frac{1}{d(\lambda_1d_1)}  \right) = \left( \frac{\partial
\rho_B}{\partial x_i} d^{2}
 +o(1) \right)\frac{|\Lambda |^2}{d^{2}} +
o\left(\frac{1}{d(\lambda_1d_1)}  \right), \, \forall\, i\in B.
$$
\noindent Thus
$$
0\geq \left( |\n \rho_B| d^{2} +o(1)   \right)\frac{|\Lambda |^2}
{d^{2}} + o\left(\frac{1}{d(\lambda_1d_1 )}  \right) \geq
C_6\frac{|\Lambda |^2}{d^{2}} + o\left(\frac{1}{d( \lambda_1d_1)}
\right) >0.
$$
This yields a contradiction. Hence \eqref{e:AB} follows.\\ The
proof of Proposition \ref{p:p1} is thereby completed.
\end{pf}

\begin{pfn}{Theorem \ref{t:11}}
Arguing by contradiction, we assume that $(P_\e)$ has a solution
whose energy is bounded. Using Theorem 1.5 of \cite{BEH} and
Proposition \ref{p:p1}, we deduce Theorem \ref{t:11}.
\end{pfn}
\section{ Appendix}
In this section, we collect some estimates needed to prove
Propositions \ref{p:214} and \ref{p:215}. Here we will denote by 
$u_\e := \sum_{j=1}^p \a_j P\d_{(x_j,\l_j)} +v_\e$ the function defined in
Theorem \ref{t:21}. Thus, we have $|x_i-x_j|\geq \a \e$ for each
$i\ne j$ and $\l_id_i \to \infty$ as $\e \to 0$ for each $i$. In
the sequel, we denote by $\var_{i,k} = \l_i^{-1}
\partial P\d_i / \partial (x_i)_k$ where $(x_i)_k$ is the $k$th
component of $x_i$, $k\in \{1,2,3\}$.\\
Recall that $B_i$ denotes $B(x_i,\a d_i /4)$ and we have, for each
$i\ne j$, $B_i\cap B_j=\emptyset$.

\begin{lem}\label{l:A0}
For $i\ne j$, we have the following estimates
\begin{align*}
1)&\quad \langle P\delta_i , \l_i \frac{\partial
P\delta_i}{\partial \l_i} \rangle=
\frac{c_1}{2}\frac{H_\varepsilon(a_i,a_i )} {\lambda_i} +
O\left(\frac{1}{ (\lambda_i d_i )^2}   \right)\\
2)&\quad\langle P\delta_j , \l_i \frac{\partial
P\delta_i}{\partial \l_i} \rangle = c_1 \left( \l_i \frac{\partial
\e_{ij}}{\partial \l_i}+ \frac{1}{2} \frac{H_\e (a_i,a_j)}{(\l_i
\l_j)^{1/2}}\right) + O\biggl(\e_{ij}^2 \left(Log
\e_{ij}^{-1}\right)^{2/3} + \sum_{k=i,j}\frac{1}
{(\l_k d_k)^2}\biggr),\\
3)&\quad\int_{A_\varepsilon}P\delta_{i}^5\lambda_i \frac{\partial
P\delta_{i} }{\partial\lambda_i} = 2\langle P\delta_{i} ,
\lambda_i\frac{\partial P\delta_{i}}{
\partial\lambda_i}\rangle + O\left(\frac{1}{(\lambda_i d_i)^2}\right),\\
4)&\quad\int_{A_\varepsilon} P\delta_j^5 \lambda_i\frac{\partial
P\delta_i}{\partial
 \lambda_i} = \langle P\delta_j , \lambda_i\frac{\partial
 P\delta_i}{\partial\lambda_i} \rangle +  O\left(\e_{ij}^2
 \left(Log \e_{ij}^{-1}\right)^{2/3} + \frac{1}
 {(\l_j d_j)^2}\right),\\
5)&\quad5\int_{A_\varepsilon}P\delta_j\left(P
\delta_i^4\lambda_i\frac{\partial P \delta_i}{\partial
\lambda_i}\right) = \langle P\delta_i , \lambda_i\frac{\partial
P\delta_i}{\partial \lambda_i}\rangle + O\left(\e_{ij}^2\left(Log
\e_{ij}^{-1} \right)^{2/3} + \frac{1}{(\l_i
d_i)^2}\right),\\
6)&\quad\int_{\R^3}\d_i ^3\d_j^3 = O\left( \e_{ij}^3
Log\e_{ij}^{-1} \right),
\end{align*}
where $c_1$ and $O$ are independent of $\e$.
\end{lem}
\begin{pf}
For the proof, we refer the interested readers to \cite{B},
\cite{R} and \cite{R2}.
\end{pf}
\begin{lem}\label{l:A1}
For $i\in \{1,...,p\}$ and $j\ne i$, we have the following
estimates
\begin{align*}
1) &\qquad  \langle P\delta_i , \var_{i,k} \rangle =
-\frac{c_1}{2\l_i ^2} \frac{\partial H_\e}{\partial (x_i)_k}
(x_i,x_i) + O\left(
\frac{1}{(\l_id_i)^3}\right),\\
2) & \qquad\int_{A_\e} P\d_i ^5 \var_{i,k} =2\langle P\delta_i ,
\var_{i,k} \rangle + O\left(
\frac{Log(\l_id_i)}{(\l_id_i)^3}\right),\\
3) & \qquad\langle P\delta_j , \var_{i,k} \rangle =
\frac{-c_1}{\l_i ^{3/2}\l_j^{1/2}} \frac{\partial H_\e}{\partial
(x_i)_k} (x_j,x_i) + \frac{c_1}{\l_i} \frac{\partial
\e_{ij}}{\partial (x_i)_k}+
O\left(\frac{1}{(\l_1d_1)^3}+\l_j |x_i-x_j|\e_{ij}^4\right),\\
4) & \qquad\int_{A_\e} P\d_j ^5 \var_{i,k} =\langle P\delta_j ,
\var_{i,k} \rangle + O\left(
\frac{1}{(\l_1d_1)^{5/2}}\right),\\
5) & \qquad\int_{A_\e} 5 P\d_j P\d_i ^4 \var_{i,k} =\langle
P\delta_j , \var_{i,k} \rangle + O\left(
\frac{1}{(\l_1d_1)^{5/2}}\right).
\end{align*}
\end{lem}
\begin{pf}
Claims 1, 2 and 3 are proved in \cite{B} and \cite{R}. We will
prove Claim 4. We have
$$
\int_{A_\e} P\d_j ^5 \var_{i,k} =\int_{A_\e} (\d_j ^5 +O(\d_j ^4
\th_j))\var_{i,k}=\langle P\delta_j , \var_{i,k} \rangle +O\left(
\int_{B_j} \d_j^4\th_j |\var_{i,k}|+\int_{A_\e \setminus B_j} \d_j
^5\d_i \right).
$$
For the second integral, using Holder's inequality, we obtain
\be\label{e:A2}
 \int_{\R^3\setminus B_j}\d_j^5\d_i= O\left(\frac{1}{(\l_jd_j)^{5/2}}\right).
\ee
 By Corollary \ref{c:A}, we have $B_i \cap B_j=\emptyset$ and
therefore, for any $x\in B_j$, we get
 \be\label{e:A3} \sup_{B_j} \biggl|
\frac{1}{\lambda_i}\frac{
\partial  \delta_i}{\partial (x_i)_k}\biggr| \leq C\sup_{B_j}
\left(\frac{1}{\l_i^{3/2}|x-x_i|^2}\right)=
O\left(\frac{1}{\l_i^{3/2}\max^2(d_i,d_j)}\right),
 \ee
 \be\label{e:A4} \sup_{B_j} \biggl| \frac{1}{\lambda_i}\frac{
\partial  \th_i}{\partial (x_i)_k}\biggr| \leq \frac{C}{\l_id_i}\sup_{B_j}
\th_i=O\left(\frac{1}{\l_i^{3/2}d_i\max(d_i,d_j)}\right). \ee
 Thus we obtain
 \be \label{e:A5}
\int_{B_j} \d_j ^4 \th_j |\var_{i,k}| \leq \frac{c}{\l_i ^{3/2}
\l_j ^{3/2} d_id_j \max(d_i,d_j) } \leq \frac{c}{(\l_1d_1)^3}.
 \ee
 Combining \eqref{e:A5} and \eqref{e:A2}, the claim follows.\\
It remains to prove Claim 5. We have
\begin{align*}
5\int_{A_\varepsilon}&P\delta_j P\delta_i^4\var_{i,k} =
5\int_{A_\e} \left(\d_i^4-4\d_i^3\th_i +
O\left(\d_i^2\th_i^2\right)\right) P\d_j
\left(\frac{1}{\lambda_i}\frac{
\partial \delta_i}{\partial (x_i)_k}-\frac{1}{\lambda_i}\frac{
\partial \th_i}{\partial (x_i)_k}\right)\\
&=\langle P\delta_j ,\var_{i,k}\rangle +O \left(\int_{B_i}\delta_j
\delta_i^4\biggl|\frac{1}{\lambda_i}\frac{
\partial \th_i}{\partial (x_i)_k}\biggr|\right)- 20 \int_{B_i}P\delta_j
\delta_i^3\th_i\frac{1}{\lambda_i}\frac{
\partial \delta_i}{\partial (x_i)_k}\\
&+ O\left(\int_{B_i}\delta_i^3\th_i^2\d_j +\int_{\R^3\setminus
B_i} \delta_i^5\d_j\right).
\end{align*}
Observe that
 \begin{align}
  \sup_{B_i}|D\th_i| & \leq
\frac{C}{d_i}\sup_{B_i}\th_i \leq \frac{C}{\l_i^{1/2}d_i^2}\quad ;
\quad \sup_{B_i}\d_j \leq \frac{c}{\l_j ^{1/2}
\max(d_i,d_j)},\label{A8}\\
\sup_{B_i}|D P\d_j| & \leq
\sup_{B_i}|D\d_j|+\sup_{B_i}|D\th_j|\leq
 \frac{C}{\l_j^{1/2} \max^2(d_i,d_j)}+
 \frac{C}{\l_j^{1/2}d_i \max(d_i,d_j)}.\label{A9}
 \end{align}
 Thus we derive
 \be\label{A5}
 \int_{B_i} \delta_i^3\th_i^2\d_j \leq |\d_j \theta_i
 ^2|_{L^\infty} \int_{B_i} \d_i ^3 \leq \frac{c
Log(\l_id_i)}{(\l_jd_j)^{1/2}(\l_i d_i)^{5/2}},
 \ee
 \be\label{A6}
 \int_{B_i}\delta_j \delta_i^4\biggl|\frac{1}{\lambda_i}\frac{
\partial \th_i}{\partial (x_i)_k}\biggr|\leq \frac{c}
{(\l_jd_j)^{1/2}(\l_i d_i)^{5/2}},
 \ee
 \be\label{A7}
 \int_{B_i}P\delta_j\delta_i^3\th_i\frac{1}{\lambda_i}\frac{
\partial \d_i}{\partial (x_i)_k} = O\left(\sup_{B_i}|D
\left(\th_iP\d_j\right)|\int_{B_i}\d_i^4|x-x_i|\right)=O\left(\frac{Log(\l_id_i)}
{(\l_id_i)^{5/2}(\l_jd_j)^{1/2}}\right).
 \ee
 Using \eqref{e:A2}, \eqref{A5}, \eqref{A6} and \eqref{A7}, the lemma follows.
\end{pf}

\begin{lem}\label{l:A2} For each $i$, we have
$$
\int_{A_\varepsilon}\biggl(\sum_{j=1}^p \a_jP\delta_j\biggr)^5
\var_{i,k}  =2\sum_{j=1}^p \a_j \langle P\delta_j ,\var_{i,k}
\rangle + O \left( \frac{1}{(\lambda_1d_1)^{9/4}}\right).
$$
\end{lem}

\begin{pf}
Notice that
\begin{align}\label{e:A6}
\biggl(\sum_{j=1}^p \a_jP\delta_j\biggr)^5&= \sum_{j=1}^p
\biggl(\a_jP\delta_j\biggr)^5
 + 5(\a_iP\d_i)^4\biggl(\sum_{j\ne i}\a_jP\d_j\biggr) +
 10(\a_iP\d_i)^3\biggl(\sum_{j\ne i}\a_jP\d_j\biggr)^2\notag\\
&+ O\biggl(\sum_{j\ne i}\d_i^2\d_j^3 + \sum_{j
\not\in\{i,r\}}\d_j^4\d_r\biggr).
\end{align}
Since $B_j\cap B_i=\emptyset$ and $B_j\cap B_r=\emptyset$, using
\eqref{e:A3} and \eqref{e:A4}, we derive
 \be\label{e:A7}
 \int_{ B_j}\d_j^4\d_r
|\var_{i,k}|\leq \frac{c}
{\l_r^{1/2}\max(d_r,d_j)\l_i^{{3/2}}d_i\max(d_i,d_j)} \int_{B_j}
\d_j ^{4}  \leq \frac{c}{(\l_1d_1)^{3}},
 \ee
 \be\label{e:A8}
 \int_{A_\e \setminus B_j}\d_j^4\d_r^2\leq \int_{A_\e \setminus
 B_j}\d_j^6 +\int_{A_\e \setminus B_j}\d_j^3\d_r^3\leq
 \frac{c}{(\l_1d_1)^{5/2}}.
 \ee
 Now we will estimate the third term. Using \eqref{A8} and
 \eqref{A9}, we obtain
 \begin{align}
\int_{B_i}P & \d_i^3\biggl(\sum_{j\ne i}\a_jP\d_j\biggr)^2
\var_{i,k}=\biggl(\sum_{j\ne i}\a_jP\d_j\biggr)^2(x_i)
\int_{B_i}\left(\d_i^3 +
O\left(\d_i^2\th_i\right)\right)\frac{1}{\l_i}\biggl(\frac{\partial
\d_i}{\partial (x_i)_k} -\frac{\partial
\th_i}{\partial(x_i)_k}\biggr)\notag\\
&+ O\biggl(\sup_{B_i}\biggl|D\biggl(\sum_{j\ne i }
\a_jP\d_j\biggr)^2\biggr| \int_{B_i}\d_i^4|x-x_i|\biggr)
=O\left(\sum_{j\ne i} \frac{Log(\l_id_i)}{\l_j\max^2(d_i,d_j)}
\frac{1}{\l_i^2 d_i} \right).\label{e:A13}
\end{align}
Combining \eqref{e:A6},...,\eqref{e:A13} and Lemma \ref{l:A1}, the
result follows.
\end{pf}

To improve the estimates of the integrals involving $v_\e$, we use
an original idea due to Rey \cite{R2}, namely we write
\be\label{A27} v_\e = \sum_{i=1}^p v_i^\e + w, \ee where $v_i^\e$
denotes the projection of $v_\e$ onto $H^1_0(B_i)$, that is
\be\label{A28} \D v_i^\e=\D v_\e\quad\mbox{in}\quad B_i;\quad
v_i^\e=0\quad\mbox{on}\quad \partial B_i, \ee where $B_i= B(x_i,\a
d_i/4)$ is defined in Corollary \ref{c:A}. $v_i^\e$ can be assumed
to be defined in all $A_\e$ since it can be continued by $0$ in
$A_\e \diagdown B_i$. We have
 \be\label{A29} v_\e =  v_i^\e + w
\quad\mbox{in}\quad B_i, \quad\mbox{with}\quad \D
w=0\,\,\mbox{in}\,\, B_i.
 \ee We split $v_i^\e$ in an even part
$v_i^{\e,e}$ and an odd part $v_i^{\e,o}$ with respect to
$(x-x_i)_k$, thus we have \be\label{A30} v_\e =
v_i^{\e,e}+v_i^{\e,o} + w \quad\mbox{in}\quad
B_i\quad\mbox{with}\quad \D w=0\,\,\mbox{in}\,\, B_i. \ee

\begin{lem}\label{l:A3}
We have
$$
\int_{B_i}\delta_i^3 v_\e^2 \frac{1}{\lambda_i}\frac{
\partial \delta_i}{\partial (x_i)_k} =O\left(||v_i^{\e,o}||
||v_\e|| + \frac{||v_\e||^2}{(\l_id_i)^{1/2}}\right).
$$

\end{lem}
\begin{pf}
Using \eqref{A30} and the fact that the even part of $v_\e^2$
has no contribution to the integrals, we obtain \be\label{A31}
\int_{B_i}\delta_i^3 v_\e^2 \frac{1}{\lambda_i}\frac{
\partial \delta_i}{\partial (x_i)_k} =\int_{B_i}\d_i^3\frac{1}{\lambda_i}\frac{
\partial \delta_i}{\partial (x_i)_k} \left(2v_\e -
w\right)w +O\left(||v_i^{\e,o}||||v_i^{\e,e}||\right).
\ee Let $\psi$ be the solution of
$$
\D\psi = \d_i^3\frac{1}{\lambda_i}\frac{
\partial \delta_i}{\partial (x_i)_k}\left(2v_\e - w\right)
\quad\mbox{in}\quad B_i;\,\,\psi=0\quad\mbox{on}\quad \partial B_i.
$$
Thus we have \be\label{A32}
\int_{B_i}\d_i^3\frac{1}{\lambda_i}\frac{
\partial \delta_i}{\partial (x_i)_k} \left(2v_\e - w\right)w=
\int_{B_i}\D\psi. w= \int_{\partial B_i}\frac{\partial\psi}{\partial\nu}w.
\ee Let $G_i$ be the Green's function for the Laplacian on $B_i$,
that is,
$$
G_i(x,y)=\frac{1}{|x-y|}-\frac{\a\d d_i}{4|x||y-\frac{(\a\d
d_i)^2x}{16|x|^2}|}, \quad (x,y)\in B_i^2.
$$
Therefore $\psi$ is given by
 \be\label{A33}
\psi(y)=-\int_{B_i}G_i(x,y)\d_i^3\frac{1}{\lambda_i}\frac{
\partial \delta_i}{\partial (x_i)_k} \left(2v_\e - w\right) dx,\quad y\in B_i
\ee and its normal derivative by
 \be\label{A34}
\frac{\partial\psi}{\partial\nu}(y)=-\int_{B_i}\frac{\partial
G_i}{\partial\nu}(x,y)\d_i^3\frac{1}{\lambda_i}\frac{
\partial \delta_i}{\partial (x_i)_k} \left(2v_\e - w\right) dx,\quad y\in \partial B_i.
 \ee
Notice that:
\begin{align}
&\mbox{for}\,\, x\in B_i\setminus B(y,{\a d_i/8}),\quad\mbox{we have}\quad \frac{\partial G_i}{\partial\nu}(x,y)= O\left(\frac{1}{d_i^2}\right),\label{A35}\\
&\mbox{for}\,\, x\in B_i\cap B(y,{\a d_i/8}),\quad\mbox{we have}\quad \frac{\partial G_i}{\partial\nu}(x,y)= O\left(\frac{1}{|x-y|^2}\right),\label{A36}\\
&\mbox{for}\,\, x\in B_i\cap B(y,{\a d_i/8}),\quad\mbox{we
have}\quad \d_i^3\frac{1}{\lambda_i}\frac{
\partial \delta_i}{\partial (x_i)_k}= O\left(\frac{1}{\l_i^2d_i^4}\right).\label{A37}
\end{align}
Therefore
\begin{align}\label{A38}
\biggl|\frac{\partial\psi}{\partial\nu}(y)\biggr|&\leq
C\int_{B_i\cap \left(|x-y|\geq {\a
d_i/8}\right)}|2v_\e-w|\frac{\d_i^4}{d_i^2}dx+
C\int_{B_i\cap\left(|x-y|\leq {\a d_i/8}\right)}
\frac{|2v_\e-w|}{\l_i^2d_i^4|x-y|^2}dx\notag\\
&\leq \frac{C}{\l_i^{1/2}d_i^2}||v_\e||,\quad \forall y\in
\partial B_i.
\end{align}
Using \eqref{A38}, \eqref{A32} becomes \be\label{A39}
\int_{B_i}\d_i^3\frac{1}{\lambda_i}\frac{
\partial \delta_i}{\partial (x_i)_k} \left(2v_\e - w\right)w=
O\left(\frac{||v_\e||}{\l_i^{1/2}d_i^2}\int_{\partial B_i}|w|\right).
\ee To estimate the right-hand side of \eqref{A39}, we introduce
the following functions
$$
\bar{w}(X)=(\a d_i/4)^{1/2}w(x_i + \a d_i X/4);\quad
\bar{v}_\e(X)= (\a d_i/4)^{1/2}v_\e(x_i + \a d_i X/4).
$$
$\bar{w}$ satisfies \be\label{A40} \D \bar{w}=0
\quad\mbox{in}\quad B:=B(0,1); \quad
\bar{w}=\bar{v}_\e\quad\mbox{on}\quad \partial B. \ee We deduce
that \be\label{A41} \int_{\partial B}|\bar{w}| \leq C\left(\int_B
|\n \bar{v}_\e|^2\right)^{1/2}=C\left(\int_{B_i} |\n
v_\e|^2\right)^{1/2}. \ee But, we have \be\label{A42}
\int_{\partial B}|\bar{w}|=\int_{\partial B}(\a d_i/4)^{1/2}|w(x_i
+ \a d_i X/4)|= \frac{1}{(\a d_i/4)^{3/2}}\int_{\partial B_i}|w|.
\ee Thus \be\label{A43} \int_{\partial B_i}|w|\leq c
d_i^{3/2}\left(\int_{B_i} |\n v_\e|^2\right)^{1/2}. \ee Using
\eqref{A31}, \eqref{A39} and \eqref{A43}, the lemma follows.
\end{pf}
\begin{lem}\label{l:A4}
For $\e$ small,  we have
$$
\int_{A_\varepsilon}\delta_i^3 v_\e \frac{1}{\lambda_i}\frac{
\partial \delta_i}{\partial (x_i)_k} =O\left(\frac{||v_i^{\e,o}||}
{\l_i^{1/2}} + \frac{||v_\e||}{\l_id_i ^{1/2}}\right).
$$
\end{lem}
\begin{pf}
Lemma \ref{l:A4} can be proved in the same way as Lemma
\ref{l:A3}. So we omit its proof.
\end{pf}
\begin{lem}\label{l:A5}
For $\e$ small and $ i\neq j,$  we have
$$
\int_{A_\varepsilon}\biggl(\sum_{j=1}^p \a_jP\delta_j\biggr)^4
v_\e \var_{i,k}=O\biggl(||v_i^{\e,o}|| \frac{1}{\l_1d_1} +
||v_\e|| \frac{1}{(\l_1d_1)^{3/2}}\biggr) .
$$
\end{lem}
\begin{pf}
We notice that
 \be\label{A15} \biggl(\sum_1^p
\a_jP\delta_j\biggr)^4=  \left(\a_iP\delta_i\right)^4
 + 4(\a_iP\d_i)^3\biggl(\sum_{j\ne i}\a_jP\d_j\biggr) +
 O\biggl(\d_i^2\sum_{j\ne i}\d_j^2 + \sum_{j \ne i}\d_j^4\biggr).
\ee
 For the last term in \eqref{A15}, we have, using \eqref{e:A2}
and \eqref{e:A3},
 \be\label{A16} \int_{A_\varepsilon}\d_j^4 v_\e
\var_{i,k} = \int_{B_j}\d_j^4 v_\e \var_{i,k} +
\int_{\R^3\setminus B_j}\d_j^4 v_\e \var_{i,k}
=O\left(\frac{||v_\e||}{\l_i^{3/2}d_i\max(d_i,d_j)\l_j^{1/2}}
+\frac{||v_\e||}{(\l_jd_j)^2}\right). \ee
 For the third term in
\eqref{A15}, we use Holder's inequality and we obtain
 \be\label{A17}
 \int_{A_\varepsilon}\d_i^2\d_j^2 |v_\e| |\var_{i,k}|\leq
 \int_{A_\varepsilon}\d_i^3\d_j^2 |v_\e| \leq c||v_\e||\e_{ij}^2
\left(Log\e_{ij}^{-1}\right)^{2/3}\leq c
\frac{||v_\e||}{(\l_1d_1)^{3/2}}.
 \ee
 Regarding the first term in \eqref{A15}, we write
 \begin{align*}
\int_{A_\varepsilon}P\d_i^4 v_\e
\var_{i,k}&=\int_{A_\e}\left(\d_i^4-4\d_i^3\th_i+
O\left(\d_i^2\th_i^2\right)\right)\frac{v_\e}{\l_i}\left(\frac{
\partial \delta_i}{\partial (x_i)_k}-\frac{
\partial \th_i}{\partial (x_i)_k}\right)\\
&=-4\th_i(x_i)\int_{B_i}\d_i^3v_\e \frac{1}{\l_i}\frac{
\partial \delta_i}{\partial (x_i)_k} + O\left(\frac{||v_\e||}{(\l_id_i)^2}\right).
\end{align*}
Using Lemma \ref{l:A4}, we derive that
 \be\label{A18}
 \int_{A_\varepsilon}P\d_i^4 v_\e \var_{i,k}=O\left(\frac{||v_i^{\e,o}||}
 {\l_id_i} + \frac{||v_\e||}{(\l_id_i)^{3/2}}\right).
 \ee
 Finally, we deal with the second term in \eqref{A15}
\begin{align}\label{A19}
\int_{B_i}P\d_i^3P\d_j v_\e
\var_{i,k}&=P\d_j(x_i)\int_{B_i}\left(\d_i^3 +
O\left(\d_i^2\th_i\right)\right)\frac{v_\e}{\l_i}\left(\frac{
\partial \delta_i}{\partial (x_i)_k}-\frac{
\partial \th_i}{\partial (x_i)_k}\right)\notag\\
&+ O\left(\sup_{B_i}|DP\d_j| \int_{B_i}\d_i^4|v_\e||x-x_i|\right).
\end{align}
Observe that, by \eqref{A8}, we have
 \be\label{A20}
P\d_j(x_i)\int_{B_i}\d_i^3\biggl(\th_i+\frac{1}{\l_i}|\frac{
\partial \th_i}{\partial (x_i)_k}|\biggr) |v_\e|=
O\left(\frac{||v_\e||}{\l_j^{1/2}d_i\max(d_i,d_j)\l_i^{3/2}}\right).
 \ee
 Using \eqref{A9}, we derive that
 \be\label{A21}
\sup_{B_i}|DP\d_j| \int_{B_i}\d_i^4|v_\e||x-x_i|=
O\left(\frac{||v_\e||}{\l_j^{1/2}d_i\max(d_i,d_j)\l_i^{3/2}}\right).
 \ee
 By Lemma \ref{l:A4}, \eqref{A20} and \eqref{A21}, \eqref{A19} become
 \be\label{A22}
\int_{B_i}P\d_i^3P\d_j v_\e \frac{1}{\lambda_i}\frac{
\partial P\delta_i}{\partial (x_i)_k}=
O\left(\frac{||v_i^{\e,o}||}{\left(\l_id_i\l_j
d_j\right)^{1/2}}+ \sum_{r\in\{i,j\}}\frac{||v_\e||}{(\l_rd_r)^{3/2}}\right).
\ee For the integral on $\R^3\setminus B_i$, we use Holder's
inequality and obtain \be\label{A23} \int_{\R^3\setminus
B_i}P\d_i^3P\d_j |v_\e| \frac{1}{\lambda_i}\frac{
\partial P\delta_i}{\partial (x_i)_k}\leq
\int_{\R^3\setminus B_i}\d_i^4\d_j |v_\e|=O\left(\frac{||v_\e||}{(\l_id_i)^2}\right).
\ee Using \eqref{A16}, \eqref{A17}, \eqref{A18}, \eqref{A22} and
\eqref{A23}, the lemma follows.
\end{pf}

\begin{lem}\label{l:210}
For  $ i\neq j$  we have
$$
\int_{A_\varepsilon}\biggl(\sum_{j=1}^p \a_jP\delta_j\biggr)^3
v_\e^2 \var_{i,k} =O\left(||v_i^{\e,o}|| ||v_\e|| +
\frac{||v_\e||^2}{(\l_1d_1)^{1/2}}\right).
$$
\end{lem}
\begin{pf}
We have
$$
\left(\sum_1^p \a_jP\delta_j\right)^3=\a_i^3\d_i^3 +
O\left(\d_i^2\th_i\right)+ O\left(\sum_{j\ne
i}\left(\d_i^2\d_j+\d_j^3\right)\right).
$$
We now observe that
 \be\label{A24} \int_{A_\e}\left(\d_i^3\th_i +
\sum_{j\ne i}\left(\d_i^3\d_j+ \d_j^3\d_i\right)\right)|v_\e|^2=
O\left(||v_\e||^2\left(\frac{1}{\l_id_i}+ \sum_{j\ne
i}\e_{ij}\left(Log\e_{ij}^{-1}\right)^{1/3}\right)\right),
 \ee
 \be\label{A25}
 \int_{\R^3\setminus B_i}\d_i^4|v_\e|^2 =
O\left(\frac{||v_\e||^2}{(\l_id_i)^2}\right), \ee \be\label{A26}
\int_{B_i}\d_i^3v_\e^2\left(\frac{1}{\lambda_i}\frac{
\partial \delta_i}{\partial (x_i)_k}-\frac{1}{\lambda_i}\frac{
\partial \th_i}{\partial (x_i)_k}\right)= O\left(||v_i^{\e,o}||||v_\e||+
\frac{||v_\e||^2}{(\l_id_i)^{1/2}}\right),
 \ee
where we have used Lemma \ref{l:A3} in the last equality. Clearly,
\eqref{A24}, \eqref{A25} and \eqref{A26} imply our lemma.
\end{pf}
\begin{lem}\label{l:213}
For $\e$ small,  we have
$$
||v_i^{\e,o}||= O \left(
\frac{1}{\left(\l_1d_1\right)^{9/8}}\right).
$$
\end{lem}
\begin{pf}
We write \be\label{A49} v_i^{\e,o}=\tilde{v}_i^o + a P\d_i +
b\l_i\frac{\partial P\d_i}{\partial\l_i} + \sum_{r=1}^3C_r
\frac{1}{\l_i}\frac{\partial P\d_i}{\partial (x_i)_r} \ee with
$$
\langle \tilde{v}_i^o, P\d_i\rangle = \langle \tilde{v}_i^o,
\frac{P\d_i}{\partial\l_i}\rangle = \langle \tilde{v}_i^o,
\frac{P\d_i}{\partial(x_i)_r}\rangle =0 \mbox{ for each } r=
1,2,3.
$$
Taking the scalar product in $H^1_0(A_\e)$ of \eqref{A49} with
$P\d_i$, $\l_i\frac{\partial P\d_i}{\partial\l_i}$,
$\frac{1}{\l_i}\frac{\partial P\d_i}{\partial(x_i)_r}$, $1\leq r
\leq 3$, provides us with the following invertible linear system
in $a$, $b$, $C_r$ (with  $1\leq r\leq 3$)
$$
(S)\,\,
\begin{cases}
\langle P\d_i, v_i^{\e,o}\rangle= a(C'+o(1))+ b\langle P\d_i,
\l_i\frac{\partial P\d_i}{\partial\l_i}\rangle +\sum_{r=1}^3C_r
\langle P\d_i, \frac{1}{\l_i}\frac{\partial P\d_i}{\partial(x_i)_r}\rangle\\
\langle\l_i\frac{\partial P\d_i}{\partial\l_i}, v_i^{\e,o}\rangle=
a \langle P\d_i, \l_i\frac{\partial P\d_i}{\partial\l_i}\rangle +
b(C'' +o(1)) + \sum_{r=1}^3C_r \langle \l_i\frac{\partial
P\d_i}{\partial\l_i},
\frac{1}{\l_i}\frac{\partial P\d_i}{\partial(x_i)_r}\rangle\\
\langle\frac{1}{\l_i}\frac{\partial P\d_i}{\partial(x_i)_j},
v_i^{\e,o}\rangle= a \langle P\d_i, \frac{1}{\l_i}\frac{\partial
P\d_i}{\partial(x_i)_j}\rangle+ b\langle \l_i\frac{\partial
P\d_i}{\partial\l_i}, \frac{1}{\l_i}\frac{\partial
P\d_i}{\partial(x_i)_j}\rangle + \sum_{r=1}^3C_r \langle
\frac{1}{\l_i}\frac{\partial P\d_i}{\partial(x_i)_j},
\frac{1}{\l_i}\frac{\partial P\d_i}{\partial(x_i)_r}\rangle.
\end{cases}
$$
Observe that
$$
\langle P\d_i, \l_i\frac{\partial P\d_i}{\partial\l_i}\rangle
=O\left(\frac{1}{\l_id_i}\right);\quad \langle \l_i\frac{\partial
P\d_i}{\partial\l_i}, \frac{1}{\l_i}\frac{\partial P\d_i}
{\partial(x_i)_r}\rangle=O\left(\frac{1}{(\l_id_i)^2}\right);
$$
$$
 \langle P\d_i, \frac{1}{\l_i}\frac{\partial
P\d_i}{\partial(x_i)_r}\rangle=
O\left(\frac{1}{(\l_id_i)^2}\right);\quad \langle
\frac{1}{\l_i}\frac{\partial P\d_i}{\partial(x_i)_j},
\frac{1}{\l_i} \frac{\partial
P\d_i}{\partial(x_i)_r}\rangle=(C'''+o(1))\d_{jr}+
O\left(\frac{1}{(\l_id_i)^2}\right),
$$
where $\d_{jr}$ denotes the Kronecker symbol.\\
 Now, because of evenness of $\d_i$ and oddness
of $v_i^{\e,o}$ with respect to $(x-x_i)_k$ we obtain
 \be\label{A50}
\langle P\d_i, v_i^{\e,o}\rangle = \int_{A_\e}\n P\d_i. \n
v_i^{\e,o} = \int_{B_i}\n P\d_i . \n v_i^{\e,o} = \int_{B_i} \d_i
^5 v_i^{\e,o} = 0.
 \ee
 In the same way we have
$$
\langle\l_i\frac{\partial P\d_i}{\partial\l_i}, v_i^{\e,o}\rangle=
\langle\frac{1}{\l_i}\frac{\partial P\d_i}{\partial(x_i)_j},
v_i^{\e,o}\rangle=0\quad\mbox{ for each } j\ne k.
$$
We also have
 \be\label{A51} \langle\var_{i,k}, v_i^{\e,o}\rangle =
\int_{B_i}\n\var_{i,k} . \n\left(v_\e-v_i^{\e,e}-w\right)
 =-\int_{A_\e\setminus B_i}\n\var_{i,k} . \n v_\e -
\int_{B_i}\n\var_{i,k} . \n w \ee
 since $v_\e$ satisfies $(V_0)$,
$v_i^{\e,e}$ is even with respect to $(x-x_i)_k$ and
$v_i^{\e,e}=0$ on $\partial B_i$. On one hand \be\label{A52}
\biggl| \int_{A_\e\setminus B_i}\n\var_{i,k}.\n v_\e\biggr| \leq C
||v_\e||\left(\int_{A_\e\setminus
B_i}|\n\var_{i,k}|^2\right)^{\frac{1}{2}}\leq
\frac{C||v_\e||}{(\l_id_i)^{\frac{3}{2}}}.
 \ee
 On the other hand, let $\psi_2$ be such that
$$
\D \psi_2=\D \var_{i,k}\, \, \mbox{in}\, \, B_i;\quad \psi_2=0\,
\, \mbox{on}\, \, \partial B_i.
$$
Writing
 \be\label{A53}
 \psi_2=  \var_{i,k} +
\th,\quad\mbox{with}\,\,\D\th=0\,\,\mbox{in}\,\,B_i,
 \ee
 we obtain
 \be\label{A54} \int_{B_i}\n\left(\var_{i,k}\right).\n w = \int_{B_i}\n\psi_2.\n w -
\int_{B_i}\n\th .\n w =-\int_{\partial
B_i}\frac{\partial\th}{\partial \nu} w.
 \ee
 Using an integral
representation for $\psi_2$, as in \eqref{A34}, we obtain for
$y\in \partial B_i$
 \be\label{A55}
 \frac{\partial\psi_2}{\partial
\nu}(y)=\int_{B_i}\frac{\partial
G_i}{\partial\nu}(x,y)\left(5\d_i^4\var_{i,k}\right)dx.
 \ee
 In $B_i\setminus B(x_i, {\a d_i}/{8})$, we argue as in
 \eqref{A38}, using \eqref{A35} and \eqref{A36}, we obtain
$$
 \int_{B_i\setminus B(x_i,{\a d_i}/{8})}\frac{\partial G_i}{\partial\nu}(x,y)
 \left(5\d_i^4\var_{i,k}\right)dx=
  O\left(\frac{1}{\l_i^{{5}/{2}}d_i^4}\right).
$$
Furthermore, since
$$
\biggl| \n \frac{\partial
G_i}{\partial\nu}(x,y)\biggr|=O\left(\frac{1}{d_i^3}\right)\quad\mbox{for}\quad
(x,y)\in B(x_i, {\a d_i/8})\times \partial B_i,
$$
we obtain
$$
\int_{ B(x_i,{\a d_i}/{8})}\frac{\partial G_i}{\partial\nu} (x,y)
\left(5\d_i^4\frac{1}{\l_i}\frac{\partial \d_i}{\partial(x_i)_k}
\right)dx \leq \frac{c}{d_i ^3}
 \int_{ B(x_i,{\a d_i}/{8})} \d_i
^5 |x-x_i|= O\left(\frac{1}{\l_i^{{3}/{2}}d_i^3}\right).
$$
where we have used the evenness of $\d_i$ and the oddness of its
derivative. Thus
 \be\label{A59} \frac{\partial\psi_2}{\partial
\nu}(y)=O\left(\frac{1}{\l_i^{3/2}d_i^3}\right)
 \ee
 so that on $\partial B_i$
 \be\label{A60}
\frac{\partial\th}{\partial\nu}=\frac{\partial\psi_2}{\partial
\nu}-\frac{\partial}{\partial\nu}\left(\frac{1}{\l_i}\frac{\partial
\d_i}{\partial(x_i)_k}\right) + \frac{\partial}{\partial\nu}
\left(\frac{1}{\l_i}\frac{\partial \th_i}{\partial(x_i)_k}\right)=
O\left(\frac{1}{\l_i^{{3}/{2}}d_i^3}\right).
 \ee
 It follows from \eqref{A51},
\eqref{A52}, \eqref{A54}, \eqref{A60} and \eqref{A43} that
\be\label{A61} \langle \frac{1}{\l_i}\frac{\partial
P\d_i}{\partial (x_i)_k}, v_i^{\e,o}\rangle =
O\left(\frac{||v_\e||}{(\l_id_i)^{3/2}}\right).
 \ee
 Inverting the
linear system $(S)$, we deduce from the above estimates
\be\label{A62}
a=O\left(\frac{||v_\e||}{(\l_id_i)^{\frac{7}{2}}}\right),\,
b=O\left(\frac{||v_\e||}{(\l_id_i)^{\frac{7}{2}}}\right), \,
C_k=O\left(\frac{||v_\e||}{(\l_id_i)^{\frac{3}{2}}}\right), \,
C_r=O\left(\frac{||v_\e||}{(\l_id_i)^{\frac{7}{2}}}\right), \,
r\ne k.
 \ee
 This implies through \eqref{A49}
 \be\label{A63}
||v_i^{\e,o}-\tilde{v}_i^o||=
O\left(\frac{||v_\e||}{(\l_id_i)^{3/2}}\right),\,\,\,
||v_i^{\e,o}||^2=||\tilde{v}_i^o||^2+
O\left(\frac{||v_\e||^2}{(\l_id_i)^3}\right).
 \ee
 We turn now to the last step, which consists in estimating
$||\tilde{v}_i^o||$. Since $\n J_\e(u_\e)=0$, we obtain
\begin{align}\label{A64}
0&=\langle \sum_{r=1}^p \a_rP\d_r + v_\e, v_i^{\e,o}\rangle -
J_\e(u_\e)^3\int_{A_\e}\left(\sum_{r=1}^p \a_rP\d_r+v_\e\right)^5 v_i^{\e,o}\\
&=\sum_{r=1}^p\a_r\int_{B_i}\d_r^5v_i^{\e,o}+ \int_{B_i}\n v_\e .
\n v_i^{\e,o}- J_\e(u_\e)^3\int_{B_i}\left(\sum_{r=1}^p
\a_rP\d_r+v_\e\right)^5 v_i^{\e,o}.\notag
\end{align}
Concerning the first integral, it is equal to $0$ if $r=i$
because of the oddness of $v_i^{\e,o}$ and the evenness of $\d_i$.
For $r\ne i$, using Holder's inequality, we obtain \be\label{A65}
\int_{B_i}\d_r^5v_i^{\e,o}=O\left(\frac{||v_i^{\e,o}||}{(\l_rd_r)^{5/2}}\right).\ee
Let us consider the second integral. Using \eqref{A30}, we obtain
\be\label{A66}
 \int_{B_i}\n v_\e . \n v_i^{\e,o}=\int_{B_i}\n
\left(v_i^{\e,o}+ v_i^{\e,e}+ w\right). \n v_i^{\e,o}=
\int_{B_i}|\n v_i^{\e,o}|^2.
 \ee
 For the last integral, we write
\begin{align}\label{A77}
\biggl(\sum_{r=1}^p \a_rP\d_r+v_\e\biggr)^5&= (\a_i P\d_i)^5 +
5(\a_iP\d_i)^4\biggl(\sum_{r\ne i} \a_rP\d_r+v_\e\biggr)\notag\\
&+ O\biggl(\d_i^3\biggl(\sum_{r\ne i}\d_r^2 + v_\e^2\biggr)+
\sum_{r\ne i}\d_r^5+ |v_\e|^5\biggr),
\end{align}
and we have to estimate the contribution of each term. We notice
that
 \be\label{A78}
 \int_{B_i}\biggl(\d_i^3\biggl(\sum_{r\ne
i}\d_r^2 + v_\e^2\biggr)+ \sum_{r\ne i}\d_r^5+
|v_\e|^5\biggr)|v_i^{\e,o}|\leq C||v_i^{\e,o}||\left( \sum
\frac{1}{(\l_jd_j)^2}+ ||v_\e||^2 \right).
 \ee
 Using \eqref{A8}, \eqref{A9} and the oddness of
$\d_i^4v_i^{\e,o}$, we obtain for $r\ne i$
\begin{align}\label{A80}
\int_{B_i}P\d_i^4P\d_r v_i^{\e,o}&=\int_{B_i}\left(\d_i^4+
O\left(\d_i^3\th_i\right)\right)\left(P\d_r(x_i)+
O\left(\frac{|x-x_i|}{\l_r^{1/2}d_i\max(d_i,d_r)}\right)\right)v_i^{\e,o}\notag\\
&=O\left(\int_{B_i}\d_i^3\th_i\d_r|v_i^{\e,o}| + \int_{B_i}\d_i^4\frac{|x-x_i||v_i^{\e,o}|}{\l_r^{1/2}d_i\max(d_i,d_r)}\right)\notag\\
&=O\left(||v_i^{\e,o}||\left(\frac{1}{(\l_id_i)^2}+
\frac{1}{(\l_rd_r)^2}\right)\right).
\end{align}
Now, we write
$$
\int_{B_i}P\d_i^4v_\e v_i^{\e,o} = \int_{B_i}P\d_i^4(v_i^{\e,o}+
v_i^{\e,e}+w) v_i^{\e,o}=\int_{B_i}P\d_i^4(v_i^{\e,o}+v_i^{\e,e})
v_i^{\e,o} + \int_{B_i}P\d_i ^4 w v_i^{\e,o}.
$$
For the first integral in the right side, we have
\begin{align}\label{A81}
\int_{B_i}P\d_i^4(v_i^{\e,o}+v_i^{\e,e})
v_i^{\e,o}&=\int_{B_i}\left(\d_i^4+O\left(\d_i^3\th_i\right)\right)
(v_i^{\e,o}+v_i^{\e,e}) v_i^{\e,o}\notag\\
&=\int_{B_i}\d_i^4(v_i^{\e,o})^2 + O\left(\frac{||v_i^{\e,o}||
||v_\e||}{\l_id_i}\right).
\end{align}
To deal with the term $\int_{B_i}P\d_i^4wv_i^{\e,o}$, we introduce
the following function
$$
\D\psi_3=P\d_i^4v_i^{\e,o}\,\,\mbox{in}\,\,B_i; \quad \psi_3=0\,\,
\mbox{on}\,\,\partial B_i.
$$
As in \eqref{A38}, we obtain
\begin{eqnarray*}
\frac{\partial\psi_3}{\partial\nu}(y)=
O\left(\frac{||v_i^{\e,o}||}{\l_i^{1/2}d_i^2}\right)\quad\mbox{for}\quad
y\in\partial B_i.
\end{eqnarray*}
Using \eqref{A43}, we find
\be\label{A83} \int_{B_i}P\d_i^4wv_i^{\e,o}= \int_{B_i} \D \psi_3
w = \int_{\partial B_i} \frac{\partial \psi_3}{\partial \nu} w
=O\left(\frac{||v_i^{\e,o}||||v_\e||}{(\l_id_i)^{1/2}}\right). \ee
Lastly, we write
\begin{align}\label{A84}
\int_{B_i}P\d_i^5 v_i^{\e,o}&= \int_{B_i}\left(\d_i^5-5\d_i^4\th_i+O\left(\d_i^3\th_i^2\right)\right)v_i^{\e,o}\notag\\
&=O\left(\sup_{B_i}|D\th_i|\int_{B_i}\d_i^4|x-x_i||v_i^{\e,o}|\right)+
 O\left(\int_{B_i}\d_i^3\th_i^2|v_i^{\e,o}|\right)
=O\left(\frac{||v_i^{\e,o}||}{(\l_id_i)^2}\right).
\end{align}
Using \eqref{A65}, ..., \eqref{A84} and the estimate of
$||v_\e||$, \eqref{A64} becomes
 \be\label{A85}
0=\int_{B_i}|\n v_i^{\e,o}|^2
-5J_\e(u_\e)^3\a_i^4\int_{B_i}\d_i^4(v_i^{\e,o})^2 +O\left(\sum
\frac{||v_i^{\e,o}||}{(\l_rd_r)^{11/8}}\right).
 \ee
Since $J_\e(u_\e)^3\a_i^4=1+o(1)$ and the quadratic form
$$
v\mapsto \int_{A_\e}|\n v|^2 - 5\int_{A_\e}\d_i^4v^2
$$
is positive definite on the subset
$\left[\mbox{Span}\,\left(P\d_i, \frac{\partial
P\d_i}{\partial\l_i}, \frac{\partial
P\d_i}{\partial(x_i)_j}\,\,1\leq j\leq
3\right)\right]^\bot_{H^1_0(A_\e)}$, we obtain
 \be\label{A88}
\int_{A_\e}|\n \tilde{v}_i^o|^2 -
5\int_{A_\e}\d_i^4(\tilde{v}_i^o)^2=O\left(\sum
\frac{1}{(\l_jd_j)^{9/4}}\right),
 \ee
 where we have used \eqref{A63}, \eqref{A85} and Proposition \ref{p:24}
and therefore our lemma follows.
\end{pf}

\end{document}